\documentclass[10pt]{article}

\title{An Explicit Shimura Tower of Function Fields over a Number Field: An Application of Takeuchi's List}
\author{
Takehiro HASEGAWA\footnote{The author was partially supported by JSPS KAKENHI Grant Number 15K17508} \footnote{thasegawa3141592@yahoo.co.jp} \\ 
Shiga University, Otsu, Shiga 520-0862, Japan \\
}

\usepackage{makeidx}
\usepackage{amsmath}
\usepackage{amssymb}
\usepackage{amsthm}
\usepackage{pxfonts}
\usepackage{enumerate}
\usepackage[all]{xy}

\makeindex

\newtheorem*{fact*}{Fact}
\newtheorem*{mainthm*}{Main theorem}

\theoremstyle{definition}
\newtheorem*{abstract*}{Abstract}
\newtheorem*{Def*}{Definition}
\newtheorem*{rem*}{Remark}
\newtheorem*{claim*}{Claim}

\newtheorem*{acknowledgement*}{Acknowledgement}

\newcommand{\Z}{{\mathbb{Z}}}
\newcommand{\Q}{{\mathbb{Q}}}

\newcommand{\C}{{\mathbb{C}}}
\newcommand{\F}{{\mathbb{F}}}
\newcommand{\Ps}{{\mathbb{P}}}
\newcommand{\As}{{\mathbb{A}}}

\begin{document}

\maketitle

\begin{abstract}
Elkies \cite{E1} proposed a procedure for constructing explicit towers of curves, and gave two towers of Shimura curves as relevant examples. 
In this paper, we present a new explicit tower of Shimura curves constructed by using this procedure. 
\end{abstract}

\bigskip 


\bigskip 


\section{Introduction}\label{intro}

According to a result obtained by Tsfasman, Vladut and Zink, excellent Goppa codes can be obtained from modular towers (more specifically, towers of modular curves) (\cite{TVZ}). 
However, the construction of such codes requires explicit modular towers. 
In 1998, Elkies introduced a procedure for defining such towers (see Proposition in \cite{E1}), and thereby constructed several classical (elliptic) modular towers, two Shimura towers and some Drinfeld modular towers (\cite{E1, E3}). 
In 2012, Hasegawa, Inuzuka and Suzuki presented a number of classical modular towers constructed by the same procedure (\cite{HIS, H}), and Garcia, Stichtenoth, Bassa and Beelen recently constructed several Drinfeld modular towers (\cite{GSBB, BB}). 
On the other hand, although Elkies suggested new Shimura towers, he did not explicitly construct any (see Section 5 in \cite{E2}). 
To the author's knowledge, only two explicit Shimura towers have been constructed thus far. 

\bigskip

Let $K$ be a totally real number field, and let $A$ be a quaternion algebra over $K$, namely, a central simple algebra over $K$ of dimension $4$. 
Throughout this paper, we assume that $A$ is ramified at all but one of the infinite places of $K$. 
Let $D(A)$ denote a discriminant of $A$, that is, the product of all finite places $\mathfrak{p}$ of $K$ such that $A$ is ramified at $\mathfrak{p}$. 
For an ideal $I$ of $K$ coprime to $D(A)$, Shimura curves $\mathcal{X}_{0}(I)$ analogous to classical modular curves $X_{0}(N)$ are well defined (see Section \ref{Elkiesexamples} below). 


\bigskip

Note that the prime number $\ell=3$ is totally ramified in $\Q(\sqrt{3})/\Q$. 
Our main theorem is as follows. 

\begin{mainthm*}
Let $K$ be the field $\Q(\sqrt{3})$, which is totally real, and let $\mathfrak{p}_{3}$ denote the place of $K$ lying above $3$. 
Then, the Shimura curve $\mathcal{X}_{0}(\mathfrak{p}_{3}^{2})$ is rational, that is, the genus of $\mathcal{X}_{0}(\mathfrak{p}_{3}^{2})$ is equal to $0$. 
For each natural number $n>1$, the curve $\mathcal{X}_{0}(\mathfrak{p}_{3}^{n})$ is defined by $n-1$ coordinates $x_{1}, x_{2}, \ldots, x_{n-1}$ which are related by the $n-2$ equations 
\begin{equation}\label{originaleq}
\left( -2 \frac{x_{j+1}-1}{x_{j+1}+2} \right)^{3}= 2 \frac{(5-3\sqrt{3})x_{j}^{3}+4}{x_{j}^{3}-2(5-3\sqrt{3})} \qquad (j=1,2, \ldots, n-2), 
\end{equation}
or, equivalently, the curve $\mathcal{X}_{0}(\mathfrak{p}_{3}^{n})$ is isomorphic to the locus of $(x_{1}, x_{2}, \ldots, x_{n-1})$ in $(\Ps^{1})^{n-1}$ satisfying the above equations  (\ref{originaleq}). 
In other words, the Shimura tower $\{\mathcal{X}_{0}(\mathfrak{p}_{3}^{n})\}_{n>1}$ is defined recursively by the affine equation 
$$
\left( -2 \frac{y-1}{y+2} \right)^{3}= 2 \frac{(5-3\sqrt{3})x^{3}+4}{x^{3}-2(5-3\sqrt{3})}. 
$$
\end{mainthm*}

\bigskip

Shimura curves $\mathcal{X}_{0}(I)$ always have Atkin-Lehner involutions $\omega_{I}$. 
Also, unlike classical modular curves, Shimura curves have no cusps. 
Thus, although the covering maps between Shimura curves cannot be computed by using $q$-expansions, they can be determined from their ramification behavior. 
In fact, Elkies illustrated this with two examples (see Third variation in \cite{E1}, and also see Section \ref{Elkiesexamples}). 
In this paper, we show our main theorem based on the method devised by Elkies (see Section \ref{original}). 

\bigskip

Beelen, Garcia and Stichtenoth derived a normal form for equations of recursive Kummer towers of function fields (\cite{BGS}). 
In their notation, the equation of our main theorem is written as 
$$
Y^{3}= \frac{-8(5-3 \sqrt{3})(X+1)^{3}+ 4(X-2)^{3}}{4(X+1)^{3}+ (5-3 \sqrt{3})(X-2)^{3}}, 
$$
where the variables $X, Y$ are defined by $X= {(2x-2)}/{(x+2)}, Y= {(2y-2)}/{(y+2)}$. 

Let $\ell \neq 2,3$ be a prime number, and let $\mathfrak{p}_{\ell}$ denote a place of $\Q(\sqrt{3})$ lying above $\ell$. 
Since the coefficients of the above equation are in the ring of integers of $\Q(\sqrt{3})$, the modulo $\mathfrak{p}_{\ell}$ reduction of our tower can be defined, and moreover, the reduction is asymptotically optimal over an extension field of $\F_{\ell}= \Z/(\ell)$. 

\bigskip

Wulftange studied one of the (modified) Shimura towers given by Elkies and proved that this tower is optimal over $\F_{2^{6}}$ (see Section 4 in \cite{W}). 
Moreover, Brander attempted to construct a Goppa code from the same tower (\cite{B}). 
In this regard, a Goppa code can be defined by using the tower in our main theorem. 

\bigskip

Let $I$ be an ideal of $K$ coprime to $D(A)$. 
For an Atkin-Lehner involution $\omega_{I}$ and a map $j \colon \mathcal{X}_{0}(I) \to \Ps_{\C}^{1}$ (e.g., see Main theorem I (3.2) in \cite{Smr}), the image of a morphism 
$$
\mathcal{X}_{0}(I) \to \Ps_{\C}^{1} \times \Ps_{\C}^{1}, \qquad z \mapsto (j(z), j(\omega_{I}(z)))
$$
is a closed subvariety of dimension $1$, and in the open affine set $(\As_{\C}^{1} \times \As_{\C}^{1}) \setminus \{\infty\}$, this variety is described by a modular polynomial $\Phi_{I}(x,y)$. 
In 1998, Elkies computed a polynomial $\Phi_{I}(x,y)$ for $N(I)=3,8,9$ (\cite{E1, E3}), where $N$ denotes the norm with respect to $K/\Q$, and furthermore, in 2005, Voight determined a polynomial $\Phi_{I}(x,y)$ for $N(I)=17$ (\cite{V1}). 
Polynomials $\Phi_{\mathfrak{p}_{3}}(x,y)$ and $\Phi_{\mathfrak{p}_{3}^{2}}(x,y)$ for our curves can be computed by using our theorem. 

\bigskip
\
This paper are organized as follows. 
In Section \ref{Elkiesexamples}, two examples of explicit Shimura towers defined by Elkies in \cite{E1} are introduced. 
In Section \ref{original}, a new explicit Shimura tower is given. 

\section{Explicit Shimura towers proposed by Elkies}\label{Elkiesexamples}

In this section, we introduce two examples by Elkies (see Third variation in \cite{E1}). 
The reason revisiting the examples of Elkies in this section is as follows: 
In \cite{E1}, Elkies left a detailed computational process for classical modular towers (see Example, First and Second variations in \cite{E1}), but he omitted the most of a computational process for Shimura towers. 
Therefore, an important relationship between Shimura towers and arithmetic groups remains unclear. 
Hence, we shall write a detailed computational process for Shimura towers, and we reveal a relationship between Shimura towers and arithmetic groups (see Remark in this section). 
Moreover, this section is also a recall of the procedure by Elkies for constructing explicit Shimura towers in the next section. 

\bigskip

Takeuchi gave a complete list of arithmetic triangle groups $\Delta$ of signature $(0; e_{1}, e_{2}, e_{3})$ (see Table (1) in \cite{T}). 
In fact, there exist 85 such groups. 
In this and the next sections, we will use three of these groups.

\bigskip

We consider the first example given by Elkies. 
Let $K$ be the number field $\Q(\sqrt{3})$, which is totally real. 
The prime numbers $\ell=3,2$ are totally ramified in $K/\Q$. 
Let $\mathfrak{p}_{3}$ (resp. $\mathfrak{p}_{2}$) denote the place of $K$ lying above $3$ (resp. $2$), that is, $\mathfrak{p}_{3}= (\sqrt{3})$ (resp. $\mathfrak{p}_{2}= (5- 3 \sqrt{3})$). 
Also, let $\Delta$ be the arithmetic triangle group of signature $(0; 2,4,12)$, and let $A$ be the quaternion algebra associated with $\Delta$, which is ramified at $\mathfrak{p}_{3}$ and at exactly one infinite place. 
Since the discriminant $D(A)$ is equal to $\mathfrak{p}_{3}$, the Shimura curves $\mathcal{X}_{0}(\mathfrak{p}_{2}^{n})$ are defined. 
Elkies subsequently constructs the explicit Shimura tower $\{\mathcal{X}_{0}(\mathfrak{p}_{2}^{n})\}_{n>1}$ as follows: 

\bigskip

Let $\mathfrak{H}$ be the upper half-plane. 
The Shimura curve $\mathcal{X}(1)= \Delta \backslash \mathfrak{H}$ is rational. 
In fact, 
$$
2g(\Delta \backslash \mathfrak{H})-2= \text{Area}(\Delta \backslash \mathfrak{H})- \sum_{j=1}^{3} \left( 1- \frac{1}{e_{J}} \right)= \frac{1}{6}- \left( \frac{1}{2}+ \frac{3}{4}+ \frac{11}{12} \right)= -2, \quad \text{and therefore} \quad g(\Delta \backslash \mathfrak{H})=0 
$$
(\cite{Smr}, and also see p.207 in \cite{T}). 
We can choose a coordinate $J$ on $\mathcal{X}(1)$ which takes the values $1, 0, \infty$ at the elliptic points $P_{2}, P_{4}, P_{12}$ of order $2,4,12$, respectively, that is, $J(P_{2})=1, J(P_{4})=0$, and $J(P_{12})=\infty$. 

In general, for each number $n$, the Shimura curve $\mathcal{X}_{0}(\mathfrak{p}_{2}^{n})$ can be identified with the set 
$$
\bigg\{ \left(P, Q \right) \in \mathcal{X}(1) \times \mathcal{X}(1) \ \bigg| \ \text{a point $P$ is $\mathfrak{p}_{2}^{n}$-isogenous to a point $Q$ } \bigg\} 
$$
(e.g., \cite{E1}), and this curve always has the Atkin-Lehner involution $\omega^{(n)} \colon (P, Q) \longleftrightarrow (Q, P)$. 
The covering map $\pi \colon \mathcal{X}_{0}(\mathfrak{p}_{2}^{n}) \to \mathcal{X}(1)$ is defined by the projection $\pi((P, Q))= P$, and its degree equals 
$$
N(\mathfrak{p}_{2}^{n}) \prod_{\mathfrak{p} \mid \mathfrak{p}_{2}^{n}} \left( 1+ \frac{1}{N(\mathfrak{p})} \right)= N(\mathfrak{p}_{2})^{n-1} (N(\mathfrak{p}_{2})+1)= 3 \cdot 2^{n-1}, 
$$
where $N$ denotes the norm with respect to $K/\Q$. 
The map $\pi$ is branched only above the elliptic points $P_{2}, P_{4}, P_{12}$ of $\mathcal{X}(1)$ and always unramified above the other points. 
The ramification index at a point $(P, Q)$ of $\mathcal{X}_{0}(\mathfrak{p}_{2}^{n})$ is equal to the denominator of the irreducible fraction 
$$
\text{ord}(Q)/\text{ord}(P), 
$$
where $\text{ord}$ represents the order of a point, and the order of a non-elliptic point is equal to 1. 

We can determine the ramification of the first map $\pi \colon \mathcal{X}_{0}(\mathfrak{p}_{2}) \to \mathcal{X}(1)$ and the involution $\omega^{(1)}$ of the curve $\mathcal{X}_{0}(\mathfrak{p}_{2})$ as follows: 
The point $P_{12}$ is totally ramified, and the point lying above $P_{12}$ is $(P_{12}, P_{4})$. 
The points lying above $P_{4}$ (resp. $P_{2}$) are $(P_{4},P_{2}), (P_{4},P_{12})$ (resp. $(P_{2},P_{4}), (P_{2},\ast)$), and thus the ramification indices are equal to $2,1$ (resp. $1,2$), respectively. 
Here, $\ast$ denotes a non-elliptic point. 
Then, the involution $\omega^{(1)}$ must interchange the points $(P_{4}, P_{12}), (P_{12}, P_{4})$ and the points $(P_{2},P_{4}), (P_{4},P_{12})$, respectively. 
Since the Hurwitz genus formula yields 
$$
2g(\mathcal{X}_{0}(\mathfrak{p}_{2}))-2= 3(-2)+ (3-1)+ (2-1)+ (2-1)=-2, \quad \text{and therefore} \quad g(\mathcal{X}_{0}(\mathfrak{p}_{2}))= 0, 
$$
the curve $\mathcal{X}_{0}(\mathfrak{p}_{2})$ is again rational, and $J$ is a polynomial of degree 3 with a triple pole, such that $J$ and $J-1$ have double zeros. 
Then, we can choose a rational coordinate $t$ on $\mathcal{X}_{0}(\mathfrak{p}_{2})$ such that 
$$
J= t(4t-3)^{2} \qquad \left(\text{so} \quad J-1= (t-1)(4t-1)^{2} \right), 
$$
and thus $\omega^{(1)}$ interchanges the points $t=0, t= \infty$ and the points $t=1, t=3/4$, respectively. 
Hence,  
\begin{equation}\label{1inv2412}
\omega^{(1)}(t)= \frac{3}{4t}. 
\end{equation} 
{\small
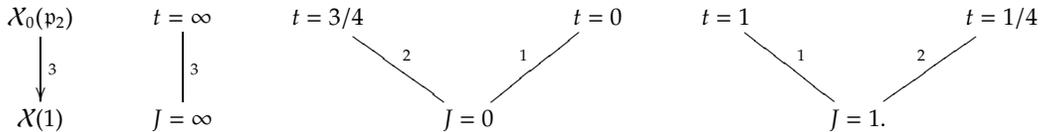
$$
\xymatrix{
\mathcal{X}_{0}(\mathfrak{p}_{2}) \ar@{->}[d]^{3} & t=\infty \ar@{-}[d]^{3} & t=3/4 \ar@{-}[dr]^{2} && t=0 \ar@{-}[dl]_{1} & t=1 \ar@{-}[dr]^{1} && t=1/4 \ar@{-}[dl]_{2} \\
\mathcal{X}(1) & J=\infty && J=0 &&& J=1. 
}
$$
}

\bigskip 

Next, we can determine the ramification of the second map $\pi \colon \mathcal{X}_{0}(\mathfrak{p}_{2}^{2}) \to \mathcal{X}_{0}(\mathfrak{p}_{2}) \to \mathcal{X}(1)$ and the involution $\omega^{(2)}$ of the curve $\mathcal{X}_{0}(\mathfrak{p}_{2}^{2})$ as follows: 
The point $P_{12}$ is totally ramified, and the point lying above $P_{12}$ is $(P_{12}, P_{2})$. 
The points lying above $P_{4}$ (resp. $P_{2}$) are $(P_{4}, \ast)$, $(P_{4},P_{4})$, $(P_{4},P_{4})$ (resp. $(P_{2},P_{12})$, $(P_{2}, P_{2})$, $(P_{2},\ast)$, $(P_{2},\ast)$), and thus the ramification indices are equal to $4,1,1$ (resp. $1,1,2,2$), respectively. 
Then, the involution $\omega^{(2)}$ interchanges the points $(P_{4}, P_{4}), (P_{4}, P_{4})$ and the points $(P_{2},P_{12}), (P_{12},P_{2})$, respectively, and fixes the point $(P_{2}, P_{2})$. 
By the Hurwitz genus formula, the curve $\mathcal{X}_{0}(\mathfrak{p}_{2}^{2})$ is rational. 
Thus, we can choose a rational coordinate $x$ on $\mathcal{X}_{0}(\mathfrak{p}_{2}^{2})$ such that 
\begin{equation}\label{eq2412}
t= \frac{x^{2}+3}{4} \qquad \left(\text{and} \quad J= \frac{x^{4}(x^{2}+3)}{4} \right), 
\end{equation}
and the involution $\omega^{(2)}$ interchanges the points $x=\sqrt{-3}, x= -\sqrt{-3}$ and the points $x=1, x= \infty$, respectively, and fixes the point $x=-1$. 
Hence,  
\begin{equation}\label{2inv2412}
\omega^{(2)}(x)= \frac{x+3}{x-1}. 
\end{equation}
{\small
$$
\xymatrix{
\mathcal{X}_{0}(\mathfrak{p}_{2}^{2}) \ar@{->}[d]^{2} & x=\infty \ar@{-}[d]^{2} & x=0 \ar@{-}[d]^{2} & x=\sqrt{-3} \ar@{-}[dr]^{1} && x=-\sqrt{-3} \ar@{-}[dl]_{1} \\
\mathcal{X}_{0}(\mathfrak{p}_{2}) & t=\infty & t=3/4 && t=0
}
$$
$$
\xymatrix{
x=1 \ar@{-}[dr]^{1} && x=-1 \ar@{-}[dl]_{1} & x=\sqrt{-2} \ar@{-}[dr]^{1} && x=-\sqrt{-2} \ar@{-}[dl]_{1} \\
& t=1 &&& t=1/4. 
}
$$
}

We now have all the information necessary to determine the curve $\mathcal{X}_{0}(\mathfrak{p}_{2}^{n})$. 
It follows from (\ref{1inv2412}) (\ref{eq2412}) (\ref{2inv2412}) that 
\begin{align*}
t \cdot \omega^{(1)}(t)= \frac{3}{4} \quad & \Leftrightarrow \quad \frac{x^{2}+3}{4} \cdot \omega^{(1)} \left( \frac{x^{2}+3}{4} \right)= \frac{3}{4} \\
\quad & \Leftrightarrow \quad (x^{2}+3) \left( \omega^{(2)}(y)^{2}+3 \right)= 12 \quad \Leftrightarrow \quad (x^{2}+3) \left( \left( \frac{y+3}{y-1} \right)^{2} +3 \right)= 12, 
\end{align*}
and hence the curve $\mathcal{X}_{0}(\mathfrak{p}_{2}^{n})$ is defined by $n-1$ coordinates $x_{1}, \ldots, x_{n-1}$ satisfying the $n-2$ relations 
$$
(x_{j}^{2}+3)(z_{j+1}^{2}+3)= 12 \qquad (j=1, \ldots, n-2), \qquad z_{j}:= (x_{j}+3)/(x_{j}-1). 
$$

Since this tower $\{ \mathcal{X}_{0}(\mathfrak{p}_{2}^{n}) \}_{n>1}$ has cyclic steps, and since it become unramified after a finite number of steps, this is dominated by a $2$-class field tower of the curve $\mathcal{X}_{0}(\mathfrak{p}_{2}^{5})$ over any finite field of odd characteristic. 

\bigskip

Next, we consider the other example given by Elkies. 
Let $K$ be the number field $\Q(\cos \frac{2}{18}  \pi)$, which is totally real, and the prime number $\ell=3$ is totally ramified in $K/\Q$. 
Also, let $\mathfrak{p}_{3}$ denote the place of $K$ lying above $3$. 
Furthermore, let $\Delta$ be the arithmetic triangle group of signature $(0; 2,3,9)$, and let $A$ be the quaternion algebra associated with $\Delta$, which is ramified at exactly two infinite places. 
Since $A$ is not ramified at all finite places of $K$, the discriminant $D(A)$ is equal to $1$, and thus the Shimura curves $\mathcal{X}_{0}(\mathfrak{p}_{3}^{n})$ are defined. 
Elkies constructs the explicit Shimura tower $\{\mathcal{X}_{0}(\mathfrak{p}_{3}^{n})\}_{n>1}$ as follows: 

\bigskip

The Shimura curve $\mathcal{X}(1)= \Delta \backslash \mathfrak{H}$ is rational. 
In fact, 
$$
2g(\Delta \backslash \mathfrak{H})-2= \text{Area}(\Delta \backslash \mathfrak{H})- \sum_{j=1}^{3} \left( 1- \frac{1}{e_{J}} \right)= \frac{1}{18}- \left( \frac{1}{2}+ \frac{2}{3}+ \frac{8}{9} \right)= -2, \quad \text{and therefore} \quad g(\Delta \backslash \mathfrak{H})=0 
$$
(\cite{Smz}, and also see p.207 in \cite{T}). 
We can choose a coordinate $J$ on $\mathcal{X}(1)$ which takes the values $1, 0, \infty$ at the elliptic points $P_{2}, P_{3}, P_{9}$ of order $2,3,9$, respectively. 
In other words, $J(P_{2})=1, J(P_{3})=0$, and $J(P_{9})=\infty$. 
For each natural number $n$, the Shimura curve $\mathcal{X}_{0}(\mathfrak{p}_{3}^{n})$ can be defined, and this curve always has the Atkin-Lehner involution $\omega^{(n)}$. 
The degree of the map $\pi \colon \mathcal{X}_{0}(\mathfrak{p}_{3}^{n}) \to \mathcal{X}(1)$ equals 
$$
N(\mathfrak{p}_{3}^{n}) \prod_{\mathfrak{p} \mid \mathfrak{p}_{3}^{n}} \left( 1+ \frac{1}{N(\mathfrak{p})} \right)= N(\mathfrak{p}_{3})^{n-1} (N(\mathfrak{p}_{3})+1)= 4 \cdot 3^{n-1}. 
$$

We can determine the ramification of the map $\mathcal{X}_{0}(\mathfrak{p}_{3}) \to \mathcal{X}(1)$ and the involution $\omega^{(1)}$ as follows: 
Note that this map is branched only above $P_{2}, P_{3}, P_{9}$ and unramified above the other points. 
The points lying above $P_{9}$ are $(P_{9},P_{9}), (P_{9},P_{3})$, and thus the ramification indices are equal to $1,3$, respectively. 
Also, the points lying above $P_{3}$ (resp. $P_{2}$) are $(P_{3},P_{9}), (P_{3},\ast)$ (resp. $(P_{2},\ast), (P_{2},\ast)$), and thus the ramification indices are equal to $1,3$  (resp. $2,2$), respectively. 
Then, the involution $\omega^{(1)}$ must interchange the points $(P_{3}, P_{9}), (P_{9}, P_{3})$ and fix the point $(P_{9},P_{9})$. 
Since, from the Hurwitz genus formula 
$$
2g(\mathcal{X}_{0}(\mathfrak{p}_{3}))-2= 4(-2)+ (3-1)+ (3-1)+ (2-1)+ (2-1)=-2, \quad \text{and therefore} \quad g(\mathcal{X}_{0}(\mathfrak{p}_{3}))= 0, 
$$
the curve $\mathcal{X}_{0}(\mathfrak{p}_{3})$ is again rational. 
Then, we can choose a rational coordinate $t$ on $\mathcal{X}_{0}(\mathfrak{p}_{3})$ such that 
$$
J= -\frac{(t-1)^{3}(9t-1)}{64 t^{3}} \qquad \left(\text{so} \quad J-1= -\frac{(3t^{2}+6t-1)^{2}}{64t^{3}} \right), 
$$
and the involution $\omega^{(1)}$ interchanges the points $t=1, t=0$ while fixing the point $t=\infty$. 
Hence,  
\begin{equation}\label{1inv239}
\omega^{(1)}(t)= 1-t. 
\end{equation} 
{\small
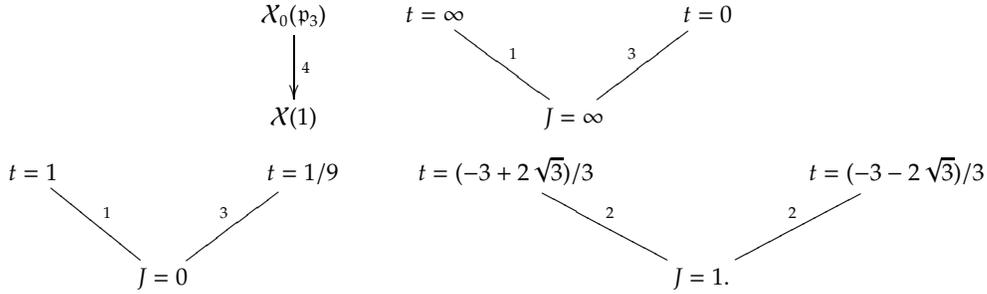
$$
\xymatrix{
\mathcal{X}_{0}(\mathfrak{p}_{3}) \ar@{->}[d]^{4} & t=\infty \ar@{-}[dr]^{1} && t=0 \ar@{-}[dl]_{3} \\
\mathcal{X}(1) && J=\infty 
}
$$
$$
\xymatrix{
t=1 \ar@{-}[dr]^{1} && t=1/9 \ar@{-}[dl]_{3} & t=(-3+2\sqrt{3})/3 \ar@{-}[dr]^{2} && t=(-3-2\sqrt{3})/3 \ar@{-}[dl]_{2} \\
& J=0 &&& J=1. 
}
$$
}

\bigskip 

Next, the ramification of the map $\mathcal{X}_{0}(\mathfrak{p}_{3}^{2}) \to \mathcal{X}_{0}(\mathfrak{p}_{3}) \to \mathcal{X}(1)$ and the involution $\omega^{(2)}$ can be determined as follows: 
The points lying above $P_{9}$ are $(P_{9}, P_{3}), (P_{9},\ast)$, and the ramification indices are $3,9$, respectively. 
The points lying above $P_{3}$ (resp. $P_{2}$) are $(P_{3}, P_{9})$, $(P_{3},P_{3})$, $(P_{3},P_{3})$, $(P_{3}, \ast)$, $(P_{3}, \ast)$, $(P_{3}, \ast)$ (resp. $(P_{2},\ast)$, $(P_{2},\ast)$, $(P_{2},\ast)$, $(P_{2},\ast)$, $(P_{2}, \ast)$, $(P_{2}, \ast)$), and therefore the ramification indices equal $1,1,1,3,3,3$ (resp. $2,2,2,2,2,2$), respectively. 
Then, the involution $\omega^{(2)}$ interchanges the points $(P_{3}, P_{9}), (P_{9}, P_{3})$ and the points $(P_{3},P_{3}), (P_{3},P_{3})$, respectively. 
By the Hurwitz genus formula, the curve $\mathcal{X}_{0}(\mathfrak{p}_{3}^{2})$ is again rational. 
Thus, we can choose a rational coordinate $x$ on $\mathcal{X}_{0}(\mathfrak{p}_{2}^{2})$ such that 
\begin{equation}\label{eq239}
t= x^{3} \qquad \left(\text{and} \quad J= -\frac{(x^{3}-1)^{3}(9x^{3}-1)}{64x^{9}} \right), 
\end{equation}
and the involution $\omega^{(2)}$ interchanges the points $x=1, x= \infty$ and the points $x= \zeta_{3}, x= \zeta_{3}^{2}$, respectively. 
Here, $\zeta_{3}$ is a cube of unity. 
Hence, the involution $\omega^{(2)}$ is given by 
\begin{equation}\label{2inv239}
\omega^{(2)}(x)= \frac{x+2}{x-1}. 
\end{equation}
{\small
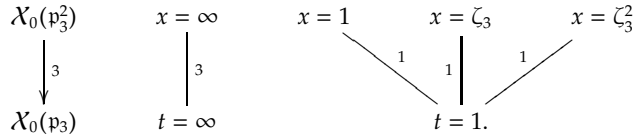
$$
\xymatrix{
\mathcal{X}_{0}(\mathfrak{p}_{3}^{2}) \ar@{->}[d]^{3} & x=\infty \ar@{-}[d]^{3} & x=1 \ar@{-}[dr]^{1} & x= \zeta_{3} \ar@{-}[d]_{1} & x= \zeta_{3}^{2} \ar@{-}[dl]_{1} \\
\mathcal{X}_{0}(\mathfrak{p}_{3}) & t=\infty && t=1. 
}
$$
}

We now have all the information necessary to determine the curve $\mathcal{X}_{0}(\mathfrak{p}_{3}^{n})$. 
It follows from (\ref{1inv239}) (\ref{eq239}) (\ref{2inv239}) that 
\begin{align*}
t+ \omega^{(1)}(t)= 1 \quad & \Leftrightarrow \quad x^{3}+ \omega^{(1)}(x^{3})= 1 \\
\quad & \Leftrightarrow \quad x^{3}+ \omega^{(2)}(y)^{3}= 1 \quad \Leftrightarrow \quad x^{3}+ \left( \frac{y+2}{y-1} \right)^{3}= 1, 
\end{align*}
and hence the curve $\mathcal{X}_{0}(\mathfrak{p}_{3}^{n})$ is defined by $n-1$ coordinates $x_{1}, \ldots, x_{n-1}$ satisfying the $n-2$ relations 
$$
x_{j}^{3}+ z_{j+1}^{3}=1 \qquad (j=1, \ldots, n-2), \qquad z_{j}:= (x_{j}+2)/(x_{j}-1). 
$$

Since this tower $\{ \mathcal{X}_{0}(\mathfrak{p}_{3}^{n}) \}_{n>1}$ has cyclic steps, and since it become unramified after finitely many steps, this is dominated by a $3$-class field tower of the curve $\mathcal{X}_{0}(\mathfrak{p}_{3}^{4})$. 

\bigskip

\begin{rem*}
It follows from Table (1) in \cite{T} that the Shimura curve $\mathcal{X}_{0}(\mathfrak{p}_{3})$ corresponds to the triangle group $\Delta$ of signature $(0; 3,3,9)$, that is, $\mathcal{X}_{0}(\mathfrak{p}_{3})= \Delta \backslash \mathfrak{H}$. 
In fact, $\Delta$ has exactly three elliptic points $(P_{9}, P_{3})$, $(P_{3}, P_{9})$, $(P_{9}, P_{9})$ of order $3,3,9$, respectively. 
Voight studied Shimura curves of genus at most $2$ (\cite{V2}). 
It follows from Table 4.3 in \cite{V2} that the curve $\mathcal{X}_{0}(\mathfrak{p}_{3}^{2})$ corresponds to a group $\Gamma$ of signature $(0;3,3,3,3)$, and the  elliptic curve $\mathcal{X}_{0}(\mathfrak{p}_{3}^{3})$ corresponds to a group of signature $(1;3,3,3)$. 
In fact, $\Gamma$ has exactly four elliptic points $(P_{3}, P_{3})$, $(P_{3}, P_{3})$, $(P_{3}, P_{9})$, $(P_{9}, P_{9})$ of order $3,3,3,3$, respectively. 
\end{rem*}

\section{A new explicit Shimura tower}\label{original}

In Section \ref{Elkiesexamples}, we introduced in detail two examples given by Elkies. 
In this section, using the method proposed by Elkies, we construct a new Shimura tower (see Main theorem in Section \ref{intro}). 

\bigskip

Let $K$ be the number field $\Q(\sqrt{3})$. 
The prime numbers $\ell=3, 2$ are totally ramified in $K/\Q$. 
Let $\mathfrak{p}_{3}$ (resp. $\mathfrak{p}_{2}$) denote the place of $K$ lying above $3$ (resp. $2$). 
Also, let $\Delta$ be the arithmetic triangle group of signature $(0; 3,3,6)$, and let $A$ be the quaternion algebra associated with $\Delta$, which is ramified at $\mathfrak{p}_{2}$ and at exactly one infinite place. 
Since the discriminant $D(A)$ equals $\mathfrak{p}_{2}$, the Shimura curves $\mathcal{X}_{0}(\mathfrak{p}_{3}^{n})$ are defined, and then the Shimura tower $\{\mathcal{X}_{0}(\mathfrak{p}_{3}^{n})\}_{n>1}$ is constructed as follows: 

\bigskip

The Shimura curve $\mathcal{X}(1)= \Delta \backslash \mathfrak{H}$ is rational. 
In fact, 
$$
2g(\Delta \backslash \mathfrak{H})-2= \text{Area}(\Delta \backslash \mathfrak{H})- \sum_{j=1}^{3} \left( 1- \frac{1}{e_{J}} \right)= \frac{1}{6}- \left( \frac{2}{3}+ \frac{2}{3}+ \frac{5}{6} \right)= -2, \quad \text{and therefore} \quad g(\Delta \backslash \mathfrak{H})=0 
$$
(\cite{Smz}, and also see p.207 in \cite{T}). 
We can choose a coordinate $J$ on $\mathcal{X}(1)$ which takes the values $1, 0, \infty$ at the elliptic points $P_{3}, P_{3}^{'}, P_{6}$ of order $3,3,6$, respectively. 
Specifically, $J(P_{3})=1, J(P_{3}^{'})=0$, and $J(P_{6})=\infty$. 
For each natural number $n$, the Shimura curve $\mathcal{X}_{0}(\mathfrak{p}_{3}^{n})$ can be defined, and this curve has the Atkin-Lehner involution $\omega^{(n)}$. 
The degree of the map $\pi \colon \mathcal{X}_{0}(\mathfrak{p}_{3}^{n}) \to \mathcal{X}(1)$ is equal to 
$$
N(\mathfrak{p}_{3}^{n}) \prod_{\mathfrak{p} \mid \mathfrak{p}_{3}^{n}} \left( 1+ \frac{1}{N(\mathfrak{p})} \right)= N(\mathfrak{p}_{3})^{n-1} (N(\mathfrak{p}_{3})+1)= 4 \cdot 3^{n-1}. 
$$

We can determine the ramification of $\mathcal{X}_{0}(\mathfrak{p}_{3}) \to \mathcal{X}(1)$ and the involution $\omega^{(1)}$ as follows: 
Note that this map is branched only above $P_{3}, P_{3}^{'}, P_{6}$ and unramified above the other points. 
The points lying above $P_{6}$ are $(P_{6},P_{3}^{'}), (P_{6},P_{3})$, and thus the ramification indices are equal to $2,2$, respectively. 
The points lying above $P_{3}^{'}$ (resp. $P_{3}$) are $(P_{3}^{'},P_{6}), (P_{3}^{'},\ast)$ (resp. $(P_{3},P_{6}), (P_{3},\ast)$), and therefore the ramification indices are equal to $1,3$  (resp. $1,3$), respectively. 
Then, the involution $\omega^{(1)}$ interchanges the points $(P_{3}^{'}, P_{6}), (P_{6}, P_{3}^{'})$ and the points $(P_{3}, P_{6}), (P_{6}, P_{3})$, respectively. 
Since by the Hurwitz genus formula 
$$
2g(\mathcal{X}_{0}(\mathfrak{p}_{3}))-2= 4(-2)+ (2-1)+ (2-1)+ (3-1)+ (3-1)=-2, \quad \text{and therefore} \quad g(\mathcal{X}_{0}(\mathfrak{p}_{3}))= 0, 
$$
the curve $\mathcal{X}_{0}(\mathfrak{p}_{3})$ is also rational. 
Then, we can choose a rational coordinate $t$ on $\mathcal{X}_{0}(\mathfrak{p}_{3})$ such that 
$$
J= \frac{4(2t+1)^{3}}{(t^{2}+10t-2)^{2}} \qquad \left(\text{so} \quad J-1= -\frac{t(t-4)^{3}}{(t^{2}+10t-2)^{2}} \right), 
$$
and the involution $\omega^{(1)}$ interchanges the points $t=\infty, t=-5+3\sqrt{3}$ and $t=0, t= -5-3\sqrt{3}$, from which 
\begin{equation}\label{1inv336}
\omega^{(1)}(t)= -\frac{(5-3\sqrt{3})t-2}{t+ (5-3\sqrt{3})}. 
\end{equation} 
{\small
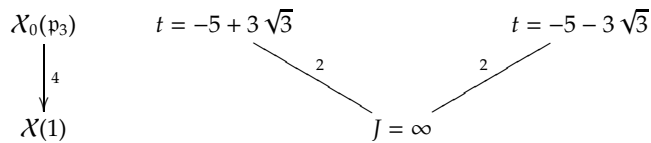
$$
\xymatrix{
\mathcal{X}_{0}(\mathfrak{p}_{3}) \ar@{->}[d]^{4} & t=-5+3\sqrt{3} \ar@{-}[dr]^{2} && t=-5-3\sqrt{3} \ar@{-}[dl]_{2} \\
\mathcal{X}(1) && J=\infty 
}
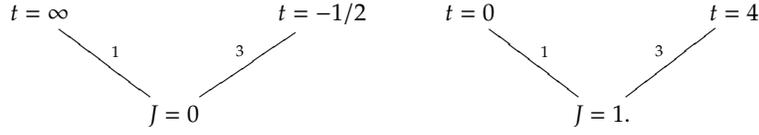
$$
$$
\xymatrix{
t=\infty \ar@{-}[dr]^{1} && t=-1/2 \ar@{-}[dl]_{3} & t=0 \ar@{-}[dr]^{1} && t=4 \ar@{-}[dl]_{3} \\
& J=0 &&& J=1. 
}
$$
}

\bigskip 

At this stage, the ramification of $\mathcal{X}_{0}(\mathfrak{p}_{3}^{2}) \to \mathcal{X}_{0}(\mathfrak{p}_{3}) \to \mathcal{X}(1)$ and the involution $\omega^{(2)}$ can be determined as follows: 
The points lying above $P_{6}$ are $(P_{6}, \ast), (P_{6}, \ast)$, and the corresponding indices are equal to $6,6$, respectively. 
The points lying above $P_{3}^{'}$ (resp. $P_{3}$) are $(P_{3}^{'}, P_{3})$, $(P_{3}^{'},P_{3}^{'})$, $(P_{3}^{'},P_{3}^{'})$, $(P_{3}^{'}, \ast)$, $(P_{3}^{'}, \ast)$, $(P_{3}^{'}, \ast)$ (resp. $(P_{3},P_{3}^{'})$, $(P_{3},P_{3})$, $(P_{3},P_{3})$, $(P_{3},\ast)$, $(P_{3}, \ast)$, $(P_{3}, \ast)$), and thus the corresponding indices are equal to $1,1,1,3,3,3$ (resp. $1,1,1,3,3,3$), respectively. 
Then, the involution $\omega^{(2)}$ interchanges the points $(P_{3}, P_{3}^{'}), (P_{3}^{'}, P_{3})$, the points $(P_{3},P_{3}), (P_{3},P_{3})$ and the points $(P_{3}^{'},P_{3}^{'}), (P_{3}^{'},P_{3}^{'})$, respectively. 
By the Hurwitz genus formula, the curve $\mathcal{X}_{0}(\mathfrak{p}_{3}^{2})$ is rational. 
Thus, we can choose a rational coordinate $x$ on $\mathcal{X}_{0}(\mathfrak{p}_{2}^{2})$ such that 
\begin{equation}\label{eq336}
t= - \frac{(5+3\sqrt{3})x^{3}+4}{x^{3}-2(5+3\sqrt{3})}, 
\end{equation}
whereby the involution $\omega^{(2)}$ interchanges the points $x=1-\sqrt{3}, x= 1+\sqrt{3}$, the points $x= (1-\sqrt{3})\zeta_{3}, x= (1-\sqrt{3})\zeta_{3}^{2}$ and the points $x= (1+\sqrt{3})\zeta_{3}, x= (1+\sqrt{3})\zeta_{3}^{2}$, respectively. 
Here, $\zeta_{3}$ is a cube of unity. 
Hence, the involution $\omega^{(2)}$ is given by 
\begin{equation}\label{2inv336}
\omega^{(2)}(x)= -2 \frac{x-1}{x+2}. 
\end{equation}
{\small
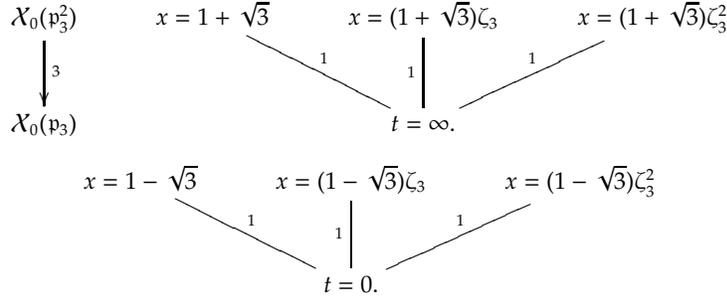
$$
\xymatrix{
\mathcal{X}_{0}(\mathfrak{p}_{3}^{2}) \ar@{->}[d]^{3} & x=1+\sqrt{3} \ar@{-}[dr]^{1} & x= (1+\sqrt{3})\zeta_{3} \ar@{-}[d]_{1} & x= (1+\sqrt{3})\zeta_{3}^{2} \ar@{-}[dl]_{1} \\
\mathcal{X}_{0}(\mathfrak{p}_{3}) && t=\infty. 
}
$$
$$
\xymatrix{
x=1-\sqrt{3} \ar@{-}[dr]^{1} & x= (1-\sqrt{3})\zeta_{3} \ar@{-}[d]_{1} & x= (1-\sqrt{3})\zeta_{3}^{2} \ar@{-}[dl]_{1} \\
& t=0. 
}
$$
}

We now have all the information necessary to determine the curve $\mathcal{X}_{0}(\mathfrak{p}_{3}^{n})$. 
It follows by (\ref{1inv336}) (\ref{eq336}) (\ref{2inv336}) that 
\begin{align*}
\omega^{(1)}(t)= & -\frac{(5-3\sqrt{3})t-2}{t+ (5-3\sqrt{3})} \quad \Leftrightarrow \quad \omega^{(1)} \left( - \frac{(5+3\sqrt{3})x^{3}+4}{x^{3}-2(5+3\sqrt{3})} \right)= \frac{4}{x^{3}} \\
\quad & \Leftrightarrow \quad - \frac{(5+3\sqrt{3}) \omega^{(2)}(y)^{3}+ 4}{\omega^{(2)}(y)^{3}- 2(5+3\sqrt{3})}= \frac{4}{x^{3}} \quad \Leftrightarrow \quad \left( -2 \frac{y-1}{y+2} \right)^{3}= 2 \frac{(5-3\sqrt{3})x^{3}+4}{x^{3}-2(5-3\sqrt{3})}, 
\end{align*}
and hence the curve $\mathcal{X}_{0}(\mathfrak{p}_{3}^{n})$ is defined by $n-1$ coordinates $x_{1}, \ldots, x_{n-1}$ satisfying the $n-2$ relations 
$$
z_{j+1}^{3}= 2 \frac{(5-3\sqrt{3})x_{j}^{3}+4}{x_{j}^{3}-2(5-3\sqrt{3})} \qquad (j=1, \ldots, n-2), \qquad z_{j}:= -2 \frac{x_{j}-1}{x_{j}+2}. 
$$

\begin{rem*}
Voight gave a list of Shimura curves of genus at most two (\cite{V2}). 
It follows from Table 4.2 in \cite{V2} that our Shimura curve $\mathcal{X}_{0}(\mathfrak{p}_{3})$ corresponds to a group $\Gamma_{1}$ of signature $(0;3,3,3,3)$, and our Shimura curve $\mathcal{X}_{0}(\mathfrak{p}_{3}^{2})$ corresponds to a group $\Gamma_{2}$ of signature $(0;3,3,3,3,3,3)$. 
In fact, $\Gamma_{1}$ has exactly four elliptic points $(P_{3},P_{6})$, $(P_{3}^{'},P_{6})$, $(P_{6},P_{3})$, $(P_{6},P_{3}^{'})$ of order $3,3,3,3$, respectively, and $\Gamma_{2}$ has exactly six elliptic points $(P_{3},P_{3})$, $(P_{3},P_{3})$, $(P_{3},P_{3}^{'})$, $(P_{3}^{'},P_{3})$, $(P_{3}^{'},P_{3}^{'})$, $(P_{3}^{'},P_{3}^{'})$ of order $3,3,3,3,3,3$, respectively. 
\end{rem*}

\begin{rem*}
While the curves given by Elkies are defined over $\Q$, ours is defined over $\Q(\sqrt{3})$. 
In general, ours also has a model over $\Q$ (see Proposition 5.1.2 in \cite{V1}). 
However, such a model does not perhaps have a simple equation as in our theorem. 
\end{rem*}

\begin{rem*}
Takeuchi gave a complete list of arithmetic triangle groups of signature $(0; e_{1}, e_{2}, e_{3})$: 
There exist 85 such groups, of which 76 groups have no cusps (see Cases II, \ldots, XIX of Table (1) in \cite{T}). 
Note that Shimura curves have no cusps. 
Elkies constructed two Shimura towers from two groups $(0; 2, 4, 12)$ and $(0; 2, 3, 9)$, respectively (\cite{E1}, and also see Section \ref{Elkiesexamples}). 
This time, we constructed a Shimura tower from the group $(0; 3, 3, 6)$ (see Section \ref{original}). 
I think that by using the same procedure as Elkies, explicit Shimura towers can not be constructed from the other 73 groups any longer. 
\end{rem*}

\end{document}